\documentclass{article}
\usepackage[round]{natbib}   

\usepackage[english]{babel}    
\usepackage{xfrac}    
\usepackage{amsfonts}
\usepackage{amsmath}
\usepackage{faktor}
\usepackage[all]{xy}
\usepackage{url}
\usepackage{xcolor}

\newtheorem{definition}{Definition}

\title{Derived symplectic geometry}
\date{}
\author{Damien Calaque \\ {\small IMAG, Univ Montpellier, CNRS, Montpellier, France} }

\begin{document}

\maketitle

\tableofcontents

\begin{abstract}
Derived symplectic geometry studies symplectic structures on derived stacks. Derived stacks are the main players in derived geometry, the purpose of which is
to deal with singular spaces, while symplectic structures are an essential ingredient of the geometric formalism of classical mechanics and classical field theory.
In addition to providing an overview of a relatively young field of research, we provide a case study on Casson's invariant. 
\end{abstract}

\noindent{\small \textbf{Keywords} -- Derived geometry. Shifted symplectic structure. Lagrangian morphism. Lagrangian intersection. Virtual class. Critical chart. Darboux theorem. Casson invariant. Donaldson--Thomas invariant. Topological field theory. }

\medskip
\colorbox[gray]{0.8}{
  \begin{minipage}{0.92\linewidth}
\textbf{Key points/objectives}
\begin{itemize}
\item Motivate the use of derived geometric techniques in symplectic geometry.  \\[-0.5cm]
\item Introduce shifted symplectic structures, with examples. \\[-0.5cm]
\item Emphasize the importance of $(-1)$-symplectic structures for virtual count. \\[-0.5cm]
\item Introduce lagrangian morphisms, and notice their relevance for symplectic geometry (moment maps, Weinstein's symplectic creed, $\dots$) and virtual count. \\[-0.5cm]
\item Provide a case study with the Casson invariant.  
\end{itemize}
  \end{minipage}}
  
\section{Introduction}

\paragraph{Derived geometry}

One of the purposes of derived geometry is to deal in a satisfying way with spaces that are singular, which is the case for many spaces appearing 
in algebraic geometry (moduli spaces) and in classical\footnote{Here, classical is understood as opposed to quantum. } physics (spaces of solutions of equations of motion). 
The reader can consult \citet{ToenDAG} and \citet{Anel} for extensive overviews of the history, results and main ideas of derived geometry. 
Foundations of derived algebraic geometry have been laid out by \citet{HAG-II} and \citet{LurieThesis}. For a treatment of (as well as references for) 
derived differential geometry we recommend the work of \citet{Pelle}. 

Very roughly, modern derived geometry combines two ways of introducing more flexibility (using homotopy theory, or higher category theory) 
into the geometry of spaces: 
\begin{itemize}
\item[(a)] It generalizes local/test objects (affine schemes, euclidean spaces, Stein manifolds, \dots) by adding homological/homotopical data to them (for instance, in algebraic 
geometry, one can replace rings with simplicial rings, or connective commutative differential graded algebras). This line of thoughts can be seen as coming out 
of an attempt to give an actual geometric content to the intersection formula from \citet{Serre}, as explained for instance by \citet{LurieThesis} and \citet{ToenDAG}. 
\item[(b)] It allows to glue local objects up to specified identifications. This led to the concept of a stack, introduced by \citet{Giraud}, and ultimately to the higher stacks of 
\citet{simpson1996algebraic}, who was somehow continuing the work of \citet{Grothendieck}. We refer to \citet{StacksSurvey} for a gentle introduction to stacks and higher stacks. 
\end{itemize}
These two ways of introducing more flexibility aim to address two different sources of singularities/pathologies: (a) bad intersections (and, more generally, 
fiber products) and (b) bad quotients. \citet{KontsevichQuasi} is one of the first authors who combined both aspects to solve enumerative problems in 
algebraic geometry, using pre-derived geometric tools. 

As a complement to the excellent introduction to the thesis of \citet{LurieThesis}, where derived (algebraic) geometry is motivated through 
the lens of Bézout's theorem (that is an intersection problem) and deformation theory, we also recommend the introduction to the paper of 
\citet{Pelle}, where the need for derived geometry is motivated by the study of solution spaces to elliptic PDEs and their local Kuranishi models. 

\paragraph{Symplectic geometry}

Symplectic geometry is a natural geometric setting for the hamiltonian formulation of classical mechanics, as one can learn for instance from 
\citet{Souriau} and \citet{Arnold}; most phase spaces appear to be symplectic manifolds (or variations of these, like Poisson manifolds). 

Recall that a symplectic manifold is a smooth manifold $X$ equipped with a $2$-form $\omega\in\Omega^2(X)$ that is non-degenerate (meaning that 
the induced bundle map $T_X\to T_X^*$ is an isomorphism) and closed (meaning that $d_{dR}\omega=0$). Observe that this definition makes sense in the 
algebro-geometric context only if $X$ is a smooth algebraic variety. The cotangent bundle $T^*M$ of a manifold $M$ is an example of a symplectic manifold, 
with $\omega_{can}=d_{dR}\lambda$, where $\lambda$ is the tautological $1$-form on $T^*M$. 
Symplectic manifolds do not have local invariants: it follows from a theorem of \cite{Darboux} that every symplectic manifold is locally symplectomorphic 
to $T^*\mathbb{R}^{n}$ equipped with $\omega_{can}$. 

Lagrangian submanifolds play a crucial role in symplectic geometry: recall that a lagrangian submanifold $L\subset X$ in a symplectic manifold $(X,\omega)$ 
is lagrangian if $\omega_{|L}=0$ and the induced map $T_L\to T^*_{L/X}$ is an isomorphism. Generalizing Darboux's theorem, \citet{WeinsteinLag} proved that in the 
neighborhood of a lagrangian submanifold $L$ every symplectic manifold is symplectomorphic to a neighborhood of the zero section of $T^*L$. 
Thus lagrangian submanifolds can naturally be interpreted as generalized configurations of a classical mechanical system. 
Lagrangian submanifolds pop up everywhere: graphs of closed $1$-forms, graphs of symplectomorphisms (an example of which is the time one flow of 
a hamiltonian vector field), conormal bundles, zero loci of moment maps, \dots This led \citet{Weinstein} to follow the symplectic creed claiming that 
``everything is a lagrangian submanifold'', and envision a symplectic category whose objects are symplectic manifolds and morphisms are lagrangian 
correspondences (i.e. lagrangian submanifolds of a product $(X_1\times X_2,p_1^*\omega_1-p_2^*\omega_2)$). 
At this point, the need to deal with singular/pathological spaces in symplectic geometry should be obvious: 
\begin{itemize}
\item The zero locus $\mu^{-1}(0)$ of a moment map $\mu$ might be singular.  
\item The above example is actually a lagrangian correspondence between the original symplectic manifold $X$ and its symplectic reduction $X_{red}$, where $X_{red}$ 
is a quotient of $\mu^{-1}(0)$ and thus could be even more singular.  
\item The composition in the symplectic category involves taking fiber products, that might not be well-behaved. 
A traditional way to deal with that is by applying a small geometric perturbation, but there are issues with this approach: 
(a) it cannot always be used in algebraic geometry and (b) geometric perturbations are not functorial. 
\end{itemize}
A leitmotiv of derived geometry is to replace geometric perturbations with homological perturbations for computing fiber products. 
Homological perturbations can be made functorial (in a higher categorical sense), and make sense in the algebro-geometric context, resolving both issues. 

\paragraph{Derived symplectic geometry}

The ancestors of shifted symplectic structures on derived stacks are the odd symplectic structures on super-manifolds and $Q$-manifolds appearing in the work of 
\citet{Schwarz}, \citet{AKSZ} and in other mathematical physics publications on the geometry of the Batalin--Vilkovisky formalism. 
These excellent works and also the beautiful treatment by \citet{Costello}, using $L_\infty$-spaces and elliptic moduli problems, have two drawbacks: 
\begin{itemize}
\item[(1)] None of the two defining properties of a symplectic structure (closedness and nondegeneracy) are homotopy invariant. 
\item[(2)] The geometric objects they consider only capture infinitesimal symmetries.  
\end{itemize}
Both issues are dealt with by the formalism of derived geometry: it indeed encompasses global symmetries (stacks have been invented for that purpose) 
and is homotopy invariant by definition. 

Let us be a bit more specific about the homotopy invariance issue: all ways of doing pre-derived geometry somehow involve spaces having a commutative differential 
graded algebra of functions, and whose tangent at a point is a cochain complex rather than just a vector space. 
Therefore, differential forms have an internal cohomological degree in addition to the usual form degree, and there are two differentials at play: the internal 
differential $\delta$ and the de Rham differential $d_{dR}$. If we have an equivalence $f:X\tilde\longrightarrow Y$ (e.g. a quasi-isomorphism of differential 
graded manifolds in the sense of \cite{Kapranov}, that is a morphism inducing an isomorphism at the level of the cohomology of the differential graded algebra 
of functions), then: 
\begin{itemize}
\item The pull-back of a $2$-form $\omega_Y$ on $Y$ that is strictly non-degenerate (meaning that the induced map $T_Y\to T^*_Y$ is an isomorphism of complexes) 
might not be strictly non-degenerate. 
\item A strictly closed $2$-form $\omega_X$ on $X$ (meaning that $d_{dR}\omega_X=\delta \omega_X=0$) may not be the pull-back of a strictly closed 
$2$-form $\omega_Y$ on $Y$.  
\end{itemize}
\citet{PTVV} introduce a very flexible notion of symplectic structure on derived stacks, that addresses all the above issues. As we see below, their 
non-degeneracy condition only requires a quasi-isomorphism, and the closedness is relaxed up to coherent homotopies (where all homotopies are 
actually part of the structure). 

In what follows, whenever we don't give a specific reference for something, the reader shall attribute it to \citet{PTVV}. 
Note that we also always work over a field of characteristic zero. 

\section{Shifted symplectic structures}

The tangent space of a derived Artin stack $X$ is no longer a vector bundle, 
but a perfect complex of $\mathcal O_X$-modules, denoted $\mathbb{T}_X$ (whose dual is denoted $\mathbb L_X$). Here $\mathcal O_X$ denotes the sheaf 
of functions on $X$. 
A $2$-form of degree $n$ on $X$ is a cochain map $\omega_0:\wedge^2\mathbb{T}_X\to\mathcal O_X[n]$. 
Note that $\omega_0$ is closed under the differential $\delta$ of the complex $\Omega^2(X)=\Gamma(X,\wedge^2\mathbb{L}_X)$ of $2$-forms on $X$; 
in particular, it induces a morphism of complexes $\mathbb{T}_X\to\mathbb{L}_X[n]$. We say that $\omega_0$ is non-degenerate if this induced morphism 
$\mathbb{T}_X\to\mathbb{L}_X[n]$ is a quasi-isomorphism, meaning that it becomes an isomorphism at the level of cohomology. 
\begin{definition}
An $n$-shifted symplectic (or, simply, $n$-symplectic) structure on $X$ is a sequence $\omega:=(\omega_0,\omega_1,\omega_2,\dots)$, where $\omega_k$ 
belongs to the complex $\Omega^k(X)$ of $k$-forms on $X$, satisfying two conditions: 
\begin{itemize}
\item The leading term $\omega_0$ is a $2$-form of degree $n$ such that is non-degenerate (as explained above, this means that the induced morphism of complexes 
$\mathbb{T}_X\to\mathbb{L}_X[n]$ becomes an isomorphism at the level of cohomology). 
\item For every $k\geq0$, $d_{dR}(\omega_k)+\delta(\omega_{k+1})=0$, meaning that $\omega_0$ is closed under the de Rham differential up to coherent 
homotopies. In other words, $\omega$ is an $n$-cocycle in the truncated 
de Rham complex 
\[
\left(\prod_{k\geq2}\Omega^k(X),d_{dR}+\delta\right)\,.
\] 
\end{itemize}
\end{definition}

\paragraph{First examples}
Because of the non-degeneracy condition, a scheme equipped with a $0$-symplectic structure must be smooth. 
Conversely, all symplectic algebraic varieties provide examples of $0$-symplectic structures. 

If $G$ is an affine algebraic group, then any non-degenerate symmetric invariant pairing $c\in S^2(\mathfrak{g}^*)^G$ on the 
Lie algebra $\mathfrak{g}$ of $G$ defines a $2$-symplectic structure on the classifying stack $BG$. 

It was proven in \citet{CalaqueCot} that the $n$-shifted cotangent stack $\mathbf{T}^*[n]X$ of an Artin stack $X$ is indeed $n$-symplectic. 
In particular, for every affine algebraic group $G$, $\mathbf{T}^*[1](BG)\simeq \big[\sfrac{\mathfrak{g}^*}{G}\big]$ is $1$-symplectic. 

Quotient stacks of quasi-symplectic groupoids, as defined by \citet{Xu}, are also $1$-symplectic, according to \citet{CalaqueSurvey}.  
In particular, for a reductive group $G$, $\big[\sfrac{G}{G^{ad}}\big]$ is $1$-symplectic. 

The moduli of objects, in the sense of \citet{ToVa}, of a $d$-Calabi--Yau category, is $(2-d)$-symplectic, 
as witnessed by \cite{ToenDAG} and shown by \cite{BD}. This implies in particular that the derived moduli 
of coherent sheaves on a Calabi--Yau $d$-fold is $(2-d)$-symplectic. 

\paragraph{The AKSZ construction}
Following \citet{AKSZ}, \citet{PTVV} developed a systematic way of producing new shifted symplectic structures on mapping stacks via transgression: 
if $(X,\omega)$ is an $n$-symplectic derived stack and $\Sigma$ is a nice enough derived stack equipped with a fundamental $d$-class $[\Sigma]$ then 
there is a $(n-d)$-symplectic structure on the mapping stack $\mathbf{Map}(\Sigma,X)$, given by $\int_{[\Sigma]}ev^*\omega$, where 
$ev:\Sigma\times \mathbf{Map}(\Sigma,X)\to X$ is the evaluation morphism. 

An important source of examples is when $\Sigma=M_B$ is the Betti stack associated with a closed oriented $d$-manifold $M$, and $G$ is a reductive group: 
then the derived stack $\mathbf{Loc}(M,G):=\mathbf{Map}(M_B,BG)$ of $G$-local systems on $M$ is $(2-d)$-symplectic\footnote{A similar result holds in derived 
differential geometry if one replaces the reductive group with a compact Lie group. }. 
Whenever $M=S^1$, $\mathbf{Loc}(S^1,G)\simeq\big[\sfrac{G}{G^{ad}}\big]$ and we recover the $1$-symplectic structure on $\big[\sfrac{G}{G^{ad}}\big]$, 
as shown by \cite{Safronov}. 

\paragraph{Virtual count on $(-1)$-symplectic derived schemes}
Let $X$ be a derived scheme equipped with a $(-1)$-shifted symplectic structure. The non-degeneracy condition for the symplectic structure imposes 
that the cotangent complex of $X$ has amplitude $[-1,0]$, meaning that $X$ is quasi-smooth. \citet{STV} proved that for a quasi-smooth derived scheme 
$X$ the cotangent complex $\mathbb{L}_X$ induces a perfect obstruction theory in the sense of \citet{BF1} on the underived truncation $t_0(X)$, 
allowing to define a well-behaved virtual fundamental class. 

\citet{PTVV} noticed that the non-degeneracy condition actually guarantees that the obstruction theory is symmetric in the sense of \cite{BF2}: in this case \citet{Behrend} 
proved that the virtual fundamental class does not depend on the choice of symmetric obstruction theory, and when $Y=t_0(X)$ is proper the virtual count (the degree 
of the virtual class) is obtained as a weighted Euler characteristic associated with a constructible function canonically defined on $Y$. 

Finally, \citet{Joyce-and-co} prove a Darboux theorem for $(-1)$-symplectic derived schemes, saying that they are locally equivalent to derived critical loci (see below). 

The above discussion for derived schemes generalizes without problem to derived Deligne--Mumford stacks (i.e. derived orbifolds).

\section{Lagrangian morphisms}

In this section we present a far reaching generalization of the notion of a lagrangian submanifold that is well-suited for derived geometry (and encompasses 
various situations): the one of a lagrangian morphism (or, more accurately, lagrangian structure on a morphism). 

Let $f:L\to X$ be a morphism of derived Artin stacks, and let $\omega$ be an $n$-symplectic structure on $X$. 
\begin{definition}
A lagrangian structure on $f$ (relatively to $\omega$) is a homotopy $\eta$ between $\omega$ and $0$  (meaning that $\omega=(d_{dR}+\delta)(\eta)$) 
in the truncated de Rham complex such that the map $\mathbb{T}_X\to\mathbb{L}_{L/X}[n-1]$ induced by $\eta_0$ is a quasi-isomorphism.  
\end{definition}
Note that the leading term $\eta_0$ of the homotopy is a homotopy between $\omega_0$ and $0$ in the complex of $2$-forms: $\omega_0=\delta\eta_0$. 
This means that the sequence $\mathbb{T}_L\to f^*\mathbb{L}_X[n]\to\mathbb{L}_L[n]$ is null-homotopic, leading to a map 
$\mathbb{T}_X\to\mathbb{L}_{L/X}[n-1]$. 

\paragraph{The derived symplectic creed}
``Everything is a lagrangian morphism''. The inclusion of a lagrangian subvariety into a symplectic variety is an example of a lagrangian morphism. 
\citet{CalaqueTFT} provides more surprising examples of lagrangian morphisms: 
\begin{itemize}
\item An $n$-symplectic structure on $X$ is the same as a lagrangian structure on $X\to *$, where $*$ is equipped with 
the trivial $(n+1)$-symplectic structure. 
\item If $G$ is an affine algebraic group, $X$ is a symplectic $G$-variety, and $\mu:X\to\mathfrak{g}^*$ is a moment map, then the induced morphism 
$\big[\sfrac{X}{G}\big]\to\big[\sfrac{\mathfrak{g}^*}{G}\big]$ between quotient stacks carries a lagrangian structure. 
\item If $G$ is reductive and $X$ is a quasi-hamiltonian $G$-variety with group-valued moment map $\mu:X\to G$, then 
$\big[\sfrac{X}{G}\big]\to\big[\sfrac{G}{G^{ad}}\big]$ also carries a lagrangian structure. 
\end{itemize}
\citet{CalaqueTFT} proves that lagrangian correspondences (i.e. lagrangian morphisms with codomain a product $(X_1\times X_2,p_1^*\omega_1-p_2^*\omega_2)$) 
compose well. This allows to define a (shifted and derived) version of the symplectic category envisioned by Weinstein, whose objects are $n$-symplectic derived stacks 
and morphisms are lagrangian correspondences. 

\paragraph{Virtual count for lagrangian intersections}
In particular, composing two lagrangian morphisms $L_1\rightarrow X\leftarrow L_2$, with $X$ being $n$-shifted symplectic, gives a lagrangian structure on 
$L_1\times_X L_2\to *$, and thus a $(n-1)$-symplectic structure on $L_1\times_X L_2$. For $n=0$ we get a $(-1)$-symplectic structure on the derived intersection 
of two lagrangian subvarieties, and the virtual count gives the intersection number of $L_1$ and $L_2$ in $X$. 
An even more specific situation is when $X=T^*M$, $L_1$ is the zero section, and $L_2$ is the graph of $d_{dR}f$ for a function $f\in \mathcal O(M)$: 
in this case the derived lagrangian intersection is the so-called derived critical locus of $f$, and the virtual count is given as the Euler characteristic of the 
hypercohomology of (a shift by $dim(M)$ of) the twisted de Rham complex $(\Omega_M,d_{dR}-d_{dR}f\wedge)$, that is a perverse sheaf on the 
reduced ordinary critical locus. \citet{SaitoSabbah} prove that, up to a shift by $dim(M)$, the twisted de Rham complex is equivalent to the sheaf of vanishing cycles 
associated with $f$, which plays a crucial role for motivic extensions and categorifications of the virtual count. 

\paragraph{The AKSZ construction as a TFT}
A relative version of the AKSZ contruction is proven by \citet{CalaqueTFT} and leads, for every $d\geq0$ and every $n$-symplectic derived stack $X$, 
to an oriented $d$-dimensional Topological Field Theory (TFT), given by $\mathbf{Map}\big((-)_B,X\big)$, and taking values in the symplectic category (with shift $n-d+1$). 
\citet{CHS} prove that this oriented TFT is fully extended in the sense of \citet{LurieTFT}. 

\section{Case study: Casson invariant}

We refer to \citet{BookCasson} for a detailed exposition, following Casson's original proposal from his talk at MSRI in 1985. 
Let $M$ be a compact oriented $3$-manifold. We consider the moduli space of irreducible $SU(2)$-representations of the fundamental group of $M$: 
\[
R^{irr}(M):=\frac{Hom^{irr}\big(\pi_1(M),SU(2)\big)}{SU(2)}\,.
\]
Irreducible representation are the ones with centralizer $\mathcal{Z}\big(SU(2)\big)=\mathbb{Z}/2\mathbb{Z}$. 
For simplicty, we assume that $M$ is an integral homology sphere (this implies in particular that the only reducible representation is the trivial one). 

\paragraph{Definition of the Casson invariant}
Let now $M=H_1\coprod_{\Sigma}H_2$ be a Heegaard splitting of $M$: $H_1$ and $H_2$ are handlebodies with boundary $\Sigma$, a closed oriented surface. 
It turns out that $R^{irr}(H_1)$ and $R^{irr}(H_2)$ are lagrangian submanifolds of $R^{irr}(\Sigma)$, that is a symplectic manifold of 
dimension $6g-6$, where $g$ is the genus of $\Sigma$. 
Even though $R^{irr}(M)$ is not a manifold, it is the intersection of $R^{irr}(H_1)$ and $R^{irr}(H_2)$, leading to define the Casson invariant $\lambda(M)$ 
as ($(-1)^g$ times) half the intersection number of $R^{irr}(H_1)$ and $R^{irr}(H_2)$ in $R^{irr}(\Sigma)$. 

\paragraph{Derived symplectic nature of the Casson invariant}
Computing the intersection number instead of looking at the naive intersection $R^{irr}(M)$ suggests one should consider the derived intersection. 
Also note that the Casson invariant is defined as half the intersection number, which is due to the hidden stacky nature of the moduli space: 
every irreducible representation carries a trivial $\mathbb{Z}/2\mathbb{Z}$ symmetry. 

The Casson invariant is in fact a virtual count on the connected component of the derived stack $\mathbf{Loc}\big(M,SU(2)\big)$ that does not contain the trivial 
local system. 
Indeed, the underived truncation of $\mathbf{Loc}\big(M,SU(2)\big)$ is $\big[\sfrac{R^{irr}(M)}{(\mathbb{Z}/2\mathbb{Z})}\big]\coprod \big[\sfrac{*}{SU(2)}\big]$. 
The $(-1)$-symplectic structure on $\mathbf{Loc}\big(M,SU(2)\big)$ therefore induces a symmetric obstruction theory on 
$\big[\sfrac{R^{irr}(M)}{(\mathbb{Z}/2\mathbb{Z})}\big]$, allowing for a virtual count. 

An important result is that $\lambda(M)$ does not depend on the choice of Heegaard splitting $M=H_1\coprod_{\Sigma}H_2$. 
This is a direct consequence of the fact that $\mathbf{Loc}\big(-,SU(2)\big)$ defines an oriented $3$-dimensional TFT with values in the derived 
version of Weinstein's symplectic category, as recalled in the previous section. 
The original definition of Casson for $\lambda(M)$ becomes a computation, that uses the fact that the 
$(-1)$-symplectic moduli can be obtained as a lagrangian intersection (thanks to the excision property of the TFT). 

 
\paragraph{Gauge theoretic approach}
\citet{Taubes} gives a gauge theoretic definition of the Casson invariant, for which we provide a derived symplectic interpretation. 
A first observation is that there is a morphism $M_{dR}\to M_{B}$ of $3$-oriented derived stacks that leads to an equivalence 
of $(-1)$-symplectic derived stacks between $\mathbf{Loc}(M,G):=\mathbf{Map}(M_B,BG)$ and $\mathbf{Flat}(M,G):=\mathbf{Map}(M_{dR},BG)$, for 
$G=SU(2)$. 

Then one identifies $\mathbf{Flat}(M,G)$ with the derived critical locus of the Chern--Simons functional $S$ that is defined on the (infinite dimensional) 
moduli stack $\mathbf{Conn}(M,G)$ of $G$-connections: 
\[
S(A):=\int_Mtr\left(d_{dR}A\wedge A+\frac23A\wedge A\wedge A\right)\,.
\]
In order to be fully accurate, one shall say that $S$ is not well-defined on this moduli space, as it is not gauge invariant: there exists a constant $c$ such that 
for a gauge $g\in C^\infty(M,G)$, 
$S(A^g)=S(g^{-1}d_{dR}g+g^{-1}Ag)=S(A)+c\mathbb{Z}$. Hence, even though $S$ is not defined on the moduli space, its exterior derivative $d_{dR}S$ 
defines a closed $1$-form on $\mathbf{Conn}(G)$, so that we can still talk about the derived critical locus of $S$ (that is the derived intersection of 
the zero section with the graph of $d_{dR}S$). 

Because of infinite dimensional issues, it is not easy at all to define correctly a (possibly virtual) count for the derived critical locus of $S$. 
This is where analytic techniques come into the game: \citet{Taubes} uses elliptic regularity and Fredholm theory to reduce locally to finite dimension. 
In fact, he constructs local potentials and finds critical charts (sometimes called Chern--Simons charts) in the sense of \cite{Joyce-and-co}, whose Darboux theorem 
in fact already provides such critical charts. 

It is \citet{tu2015casson} who first gave a virtual count formulation of the gauge theoretic approach to the Casson invariant, using $L_\infty$-spaces and Kuranishi charts. 

\section{Conclusion}

An important topic of current research is the quantization of shifted symplectic structures. The study of their geometric quantization has been initiated by 
\citet{safronov2020shifted}, while the one of their deformation quantization has been initiated by \citet{CPTVV}. The latter required to develop a whole new 
theory of shifted Poisson structure, which was independently done by \citet{Pridham}. 

The virtual count on $(-1)$-symplectic derived schemes can already be seen as a quantization process. General quantization principles suggest that a 
$(-1)$-symplectic structures shall be quantized by a cochain complex, whose Euler characteristic would give back the virtual count. 
\citet{Joyce-perverse} prove that if a scheme $X_0$ admits local critical charts with a global choice of appropriate signs for their glueing (such an $X_0$ is called 
an oriented d-critical locus), 
then one can construct a perverse sheaf on $X_0$ that is locally equivalent to the sheaf of vanishing cycles of the local potential; the Euler characteristic of the 
hypercohomology of this sheaf gives back the virtual count on $X_0$. 

Examples of such oriented d-critical loci are given as underived truncations of $(-1)$-symplectic derived schemes equipped with a square root of the canonical sheaf 
(such a square root exists for instance whenever there is a lagrangian foliation structure on the $(-1)$-symplectic derived scheme). \citet{Pridham-Q} provides a 
detailed explanation of the relation between the above and deformation quantization of $(-1)$-symplectic structures. 

This circle of ideas can be applied to recover a holomorphic version of the Casson invariant, originally introduced by \citet{Thomas}, that ``counts'' sheaves 
on a Calabi--Yau $3$-fold, and that is the starting point of the theory of Donaldson--Thomas invariants.

\paragraph{Acknowledgments}

This project has received funding from the European
Research Council (ERC) under the European Union’s Horizon 2020 research and innovation
programme (Grant Agreement No. 768679).

\bibliographystyle{apa-good}

\bibliography{my_bibtex}

\begin{thebibliography}{45}
\expandafter\ifx\csname natexlab\endcsname\relax\def\natexlab#1{#1}\fi
\expandafter\ifx\csname url\endcsname\relax
  \def\url#1{{\tt #1}}\fi
\expandafter\ifx\csname urlprefix\endcsname\relax\def\urlprefix{URL }\fi

\bibitem[{Akbulut \& McCarthy(1990)}]{BookCasson}
Akbulut, S., \& McCarthy, J.~D. (1990).
\newblock {\em Casson's invariant for oriented homology {$3$}-spheres\/},
  vol.~36 of {\em Mathematical Notes\/}.
\newblock Princeton University Press, Princeton, NJ.
\newline\urlprefix\url{https://doi.org/10.1515/9781400860623}

\bibitem[{Alexandrov et~al.(1997)Alexandrov, Schwarz, Zaboronsky, \&
  Kontsevich}]{AKSZ}
Alexandrov, M., Schwarz, A., Zaboronsky, O., \& Kontsevich, M. (1997).
\newblock The geometry of the master equation and topological quantum field
  theory.
\newblock {\em Internat. J. Modern Phys. A\/}, {\em 12\/}(7), 1405--1429.
\newline\urlprefix\url{https://doi.org/10.1142/S0217751X97001031}

\bibitem[{Anel(2021)}]{Anel}
Anel, M. (2021).
\newblock The geometry of ambiguity: an introduction to the ideas of derived
  geometry.
\newblock In {\em New spaces in mathematics---formal and conceptual
  reflections\/}, (pp. 505--553). Cambridge Univ. Press, Cambridge.

\bibitem[{Arnol'd(1989)}]{Arnold}
Arnol'd, V.~I. (1989).
\newblock {\em Mathematical methods of classical mechanics\/}, vol.~60 of {\em
  Graduate Texts in Mathematics\/}.
\newblock Springer-Verlag, New York, second ed.
\newblock Translated from the Russian by K. Vogtmann and A. Weinstein.
\newline\urlprefix\url{https://doi.org/10.1007/978-1-4757-2063-1}

\bibitem[{Behrend(2009)}]{Behrend}
Behrend, K. (2009).
\newblock Donaldson-{T}homas type invariants via microlocal geometry.
\newblock {\em Ann. of Math. (2)\/}, {\em 170\/}(3), 1307--1338.
\newline\urlprefix\url{https://doi.org/10.4007/annals.2009.170.1307}

\bibitem[{Behrend \& Fantechi(1997)}]{BF1}
Behrend, K., \& Fantechi, B. (1997).
\newblock The intrinsic normal cone.
\newblock {\em Invent. Math.\/}, {\em 128\/}(1), 45--88.
\newline\urlprefix\url{https://doi.org/10.1007/s002220050136}

\bibitem[{Behrend \& Fantechi(2008)}]{BF2}
Behrend, K., \& Fantechi, B. (2008).
\newblock Symmetric obstruction theories and {H}ilbert schemes of points on
  threefolds.
\newblock {\em Algebra Number Theory\/}, {\em 2\/}(3), 313--345.
\newline\urlprefix\url{https://doi.org/10.2140/ant.2008.2.313}

\bibitem[{Brav et~al.(2015)Brav, Bussi, Dupont, Joyce, \&
  Szendr\H{o}i}]{Joyce-perverse}
Brav, C., Bussi, V., Dupont, D., Joyce, D., \& Szendr\H{o}i, B. (2015).
\newblock Symmetries and stabilization for sheaves of vanishing cycles.
\newblock {\em J. Singul.\/}, {\em 11\/}, 85--151.
\newblock With an appendix by J\"{o}rg Sch\"{u}rmann.
\newline\urlprefix\url{https://doi.org/10.5427/jsing.2015.11e}

\bibitem[{Brav et~al.(2019)Brav, Bussi, \& Joyce}]{Joyce-and-co}
Brav, C., Bussi, V., \& Joyce, D. (2019).
\newblock A {D}arboux theorem for derived schemes with shifted symplectic
  structure.
\newblock {\em J. Amer. Math. Soc.\/}, {\em 32\/}(2), 399--443.
\newline\urlprefix\url{https://doi.org/10.1090/jams/910}

\bibitem[{Brav \& Dyckerhoff(2021)}]{BD}
Brav, C., \& Dyckerhoff, T. (2021).
\newblock Relative {C}alabi-{Y}au structures {II}: shifted {L}agrangians in the
  moduli of objects.
\newblock {\em Selecta Math. (N.S.)\/}, {\em 27\/}(4), Paper No. 63, 45.
\newline\urlprefix\url{https://doi.org/10.1007/s00029-021-00642-5}

\bibitem[{Calaque(2015)}]{CalaqueTFT}
Calaque, D. (2015).
\newblock Lagrangian structures on mapping stacks and semi-classical {TFT}s.
\newblock In {\em Stacks and categories in geometry, topology, and algebra\/},
  vol. 643 of {\em Contemp. Math.\/}, (pp. 1--23). Amer. Math. Soc.,
  Providence, RI.
\newline\urlprefix\url{https://doi.org/10.1090/conm/643/12894}

\bibitem[{Calaque(2019)}]{CalaqueCot}
Calaque, D. (2019).
\newblock Shifted cotangent stacks are shifted symplectic.
\newblock {\em Ann. Fac. Sci. Toulouse Math. (6)\/}, {\em 28\/}(1), 67--90.
\newline\urlprefix\url{https://doi.org/10.5802/afst.1593}

\bibitem[{Calaque(2021)}]{CalaqueSurvey}
Calaque, D. (2021).
\newblock Derived stacks in symplectic geometry.
\newblock In {\em New spaces in physics---formal and conceptual reflections\/},
  (pp. 155--201). Cambridge Univ. Press, Cambridge.

\bibitem[{Calaque et~al.(2022)Calaque, Haugseng, \& Scheimbauer}]{CHS}
Calaque, D., Haugseng, R., \& Scheimbauer, C. (2022).
\newblock The {AKSZ} construction in derived algebraic geometry as an extended
  topological field theory.
\newline\urlprefix\url{https://doi.org/10.48550/arXiv.2108.02473}

\bibitem[{Calaque et~al.(2017)Calaque, Pantev, To\"{e}n, Vaqui\'{e}, \&
  Vezzosi}]{CPTVV}
Calaque, D., Pantev, T., To\"{e}n, B., Vaqui\'{e}, M., \& Vezzosi, G. (2017).
\newblock Shifted {P}oisson structures and deformation quantization.
\newblock {\em J. Topol.\/}, {\em 10\/}(2), 483--584.
\newline\urlprefix\url{https://doi.org/10.1112/topo.12012}

\bibitem[{Costello(2013)}]{Costello}
Costello, K. (2013).
\newblock Notes on supersymmetric and holomorphic field theories in dimensions
  2 and 4.
\newblock {\em Pure Appl. Math. Q.\/}, {\em 9\/}(1), 73--165.
\newline\urlprefix\url{https://doi.org/10.4310/PAMQ.2013.v9.n1.a3}

\bibitem[{Darboux(1882)}]{Darboux}
Darboux, G. (1882).
\newblock Sur le problème de {P}faff.
\newblock {\em Bulletin des Sciences Mathématiques et Astronomiques, Série
  2\/}, {\em 6\/}(1), 14--36.
\newline\urlprefix\url{http://www.numdam.org/item/BSMA_1882_2_6_1_14_1/}

\bibitem[{Giraud(1971)}]{Giraud}
Giraud, J. (1971).
\newblock {\em Cohomologie non ab\'{e}lienne\/}.
\newblock Die Grundlehren der mathematischen Wissenschaften, Band 179.
  Springer-Verlag, Berlin-New York.

\bibitem[{Grothendieck(1983)}]{Grothendieck}
Grothendieck, A. (1983).
\newblock {\em Pursuing stacks (\`a la poursuite des champs). {V}ol. {I}\/},
  vol.~20 of {\em Documents Math\'{e}matiques (Paris) [Mathematical Documents
  (Paris)]\/}.
\newblock Soci\'{e}t\'{e} Math\'{e}matique de France, Paris, \copyright 2022.
\newblock Edited by Georges Maltsiniotis.

\bibitem[{Kapranov(2001)}]{Kapranov}
Kapranov, M. (2001).
\newblock Injective resolutions of {$BG$} and derived moduli spaces of local
  systems.
\newblock {\em J. Pure Appl. Algebra\/}, {\em 155\/}(2-3), 167--179.
\newline\urlprefix\url{https://doi.org/10.1016/S0022-4049(99)00109-7}

\bibitem[{Kontsevich(1995)}]{KontsevichQuasi}
Kontsevich, M. (1995).
\newblock Enumeration of rational curves via torus actions.
\newblock In {\em The moduli space of curves ({T}exel {I}sland, 1994)\/}, vol.
  129 of {\em Progr. Math.\/}, (pp. 335--368). Birkh\"{a}user Boston, Boston,
  MA.
\newline\urlprefix\url{https://doi.org/10.1007/978-1-4612-4264-2_12}

\bibitem[{Lurie(2004)}]{LurieThesis}
Lurie, J. (2004).
\newblock {\em Derived algebraic geometry\/}.
\newblock ProQuest LLC, Ann Arbor, MI.
\newblock Thesis (Ph.D.)--Massachusetts Institute of Technology.
\newline\urlprefix\url{http://hdl.handle.net/1721.1/30144}

\bibitem[{Lurie(2009)}]{LurieTFT}
Lurie, J. (2009).
\newblock On the classification of topological field theories.
\newblock In {\em Current developments in mathematics, 2008\/}, (pp. 129--280).
  Int. Press, Somerville, MA.

\bibitem[{Mestrano \& Simpson(2021)}]{StacksSurvey}
Mestrano, N., \& Simpson, C. (2021).
\newblock Stacks.
\newblock In {\em New spaces in mathematics---formal and conceptual
  reflections\/}, (pp. 462--504). Cambridge Univ. Press, Cambridge.

\bibitem[{Pantev et~al.(2013)Pantev, To\"{e}n, Vaqui\'{e}, \& Vezzosi}]{PTVV}
Pantev, T., To\"{e}n, B., Vaqui\'{e}, M., \& Vezzosi, G. (2013).
\newblock Shifted symplectic structures.
\newblock {\em Publ. Math. Inst. Hautes \'{E}tudes Sci.\/}, {\em 117\/},
  271--328.
\newline\urlprefix\url{https://doi.org/10.1007/s10240-013-0054-1}

\bibitem[{Pridham(2017)}]{Pridham}
Pridham, J.~P. (2017).
\newblock Shifted {P}oisson and symplectic structures on derived {$N$}-stacks.
\newblock {\em J. Topol.\/}, {\em 10\/}(1), 178--210.
\newline\urlprefix\url{https://doi.org/10.1112/topo.12004}

\bibitem[{Pridham(2019)}]{Pridham-Q}
Pridham, J.~P. (2019).
\newblock Deformation quantisation for {$(-1)$}-shifted symplectic structures
  and vanishing cycles.
\newblock {\em Algebr. Geom.\/}, {\em 6\/}(6), 747--779.
\newline\urlprefix\url{https://doi.org/10.14231/ag-2019-032}

\bibitem[{Sabbah \& Saito(2014)}]{SaitoSabbah}
Sabbah, C., \& Saito, M. (2014).
\newblock Kontsevich's conjecture on an algebraic formula for vanishing cycles
  of local systems.
\newblock {\em Algebr. Geom.\/}, {\em 1\/}(1), 107--130.
\newline\urlprefix\url{https://doi.org/10.14231/AG-2014-006}

\bibitem[{Safronov(2016)}]{Safronov}
Safronov, P. (2016).
\newblock Quasi-{H}amiltonian reduction via classical {C}hern-{S}imons theory.
\newblock {\em Adv. Math.\/}, {\em 287\/}, 733--773.
\newline\urlprefix\url{https://doi.org/10.1016/j.aim.2015.09.031}

\bibitem[{Safronov(2020)}]{safronov2020shifted}
Safronov, P. (2020).
\newblock Shifted geometric quantization.
\newline\urlprefix\url{https://doi.org/10.48550/arXiv.2011.05730}

\bibitem[{Sch\"{u}rg et~al.(2015)Sch\"{u}rg, To\"{e}n, \& Vezzosi}]{STV}
Sch\"{u}rg, T., To\"{e}n, B., \& Vezzosi, G. (2015).
\newblock Derived algebraic geometry, determinants of perfect complexes, and
  applications to obstruction theories for maps and complexes.
\newblock {\em J. Reine Angew. Math.\/}, {\em 702\/}, 1--40.
\newline\urlprefix\url{https://doi.org/10.1515/crelle-2013-0037}

\bibitem[{Schwarz(1993)}]{Schwarz}
Schwarz, A. (1993).
\newblock Geometry of {B}atalin-{V}ilkovisky quantization.
\newblock {\em Comm. Math. Phys.\/}, {\em 155\/}(2), 249--260.
\newline\urlprefix\url{http://projecteuclid.org/euclid.cmp/1104253279}

\bibitem[{Serre(1965)}]{Serre}
Serre, J.-P. (1965).
\newblock {\em Alg\`ebre locale. {M}ultiplicit\'{e}s\/}, vol.~11 of {\em
  Lecture Notes in Mathematics\/}.
\newblock Springer-Verlag, Berlin-New York.
\newblock Cours au Coll\`ege de France, 1957--1958, r\'{e}dig\'{e} par Pierre
  Gabriel, Seconde \'{e}dition, 1965.

\bibitem[{Simpson(1996)}]{simpson1996algebraic}
Simpson, C. (1996).
\newblock Algebraic (geometric) $n$-stacks.
\newline\urlprefix\url{https://doi.org/10.48550/arXiv.alg-geom/9609014}

\bibitem[{Souriau(1970)}]{Souriau}
Souriau, J.-M. (1970).
\newblock {\em Structure des syst\`emes dynamiques\/}.
\newblock Dunod, Paris.
\newblock Ma\^{i}trises de math\'{e}matiques.

\bibitem[{Steffens(2023)}]{Pelle}
Steffens, P. (2023).
\newblock Derived {$C^{\infty}$}-geometry {I}: {F}oundations.
\newline\urlprefix\url{https://doi.org/10.48550/arXiv.2304.08671}

\bibitem[{Taubes(1990)}]{Taubes}
Taubes, C.~H. (1990).
\newblock Casson's invariant and gauge theory.
\newblock {\em J. Differential Geom.\/}, {\em 31\/}(2), 547--599.
\newline\urlprefix\url{http://projecteuclid.org/euclid.jdg/1214444327}

\bibitem[{Thomas(2000)}]{Thomas}
Thomas, R.~P. (2000).
\newblock A holomorphic {C}asson invariant for {C}alabi-{Y}au 3-folds, and
  bundles on {$K3$} fibrations.
\newblock {\em J. Differential Geom.\/}, {\em 54\/}(2), 367--438.
\newline\urlprefix\url{http://projecteuclid.org/euclid.jdg/1214341649}

\bibitem[{To\"{e}n(2014)}]{ToenDAG}
To\"{e}n, B. (2014).
\newblock Derived algebraic geometry.
\newblock {\em EMS Surv. Math. Sci.\/}, {\em 1\/}(2), 153--240.
\newline\urlprefix\url{https://doi.org/10.4171/EMSS/4}

\bibitem[{To\"{e}n \& Vaqui\'{e}(2007)}]{ToVa}
To\"{e}n, B., \& Vaqui\'{e}, M. (2007).
\newblock Moduli of objects in dg-categories.
\newblock {\em Ann. Sci. \'{E}cole Norm. Sup. (4)\/}, {\em 40\/}(3), 387--444.
\newline\urlprefix\url{https://doi.org/10.1016/j.ansens.2007.05.001}

\bibitem[{To\"{e}n \& Vezzosi(2008)}]{HAG-II}
To\"{e}n, B., \& Vezzosi, G. (2008).
\newblock Homotopical algebraic geometry. {II}. {G}eometric stacks and
  applications.
\newblock {\em Mem. Amer. Math. Soc.\/}, {\em 193\/}(902), x+224.
\newline\urlprefix\url{https://doi.org/10.1090/memo/0902}

\bibitem[{Tu(2015)}]{tu2015casson}
Tu, J. (2015).
\newblock Casson invariants via virtual counting.
\newline\urlprefix\url{https://doi.org/10.48550/arXiv.1512.02322}

\bibitem[{Weinstein(1971)}]{WeinsteinLag}
Weinstein, A. (1971).
\newblock Symplectic manifolds and their {L}agrangian submanifolds.
\newblock {\em Advances in Math.\/}, {\em 6\/}, 329--346 (1971).
\newline\urlprefix\url{https://doi.org/10.1016/0001-8708(71)90020-X}

\bibitem[{Weinstein(1981)}]{Weinstein}
Weinstein, A. (1981).
\newblock Symplectic geometry.
\newblock {\em Bull. Amer. Math. Soc. (N.S.)\/}, {\em 5\/}(1), 1--13.
\newline\urlprefix\url{https://doi.org/10.1090/S0273-0979-1981-14911-9}

\bibitem[{Xu(2004)}]{Xu}
Xu, P. (2004).
\newblock Momentum maps and {M}orita equivalence.
\newblock {\em J. Differential Geom.\/}, {\em 67\/}(2), 289--333.
\newline\urlprefix\url{http://projecteuclid.org/euclid.jdg/1102536203}

\end{thebibliography}
\end{document}